\documentclass[paper=a4,fontsize=12pt,abstract]{scrartcl}
\pdfoutput=1 
\usepackage[ngerman,english]{babel,varioref}
\usepackage[utf8]{inputenc} 
\usepackage[T1]{fontenc} 
\usepackage{lmodern,microtype} 
\usepackage{csquotes} 

\usepackage{amssymb}
\usepackage{mathtools}
\usepackage{amsthm}
\usepackage{pgf,tikz}
\usetikzlibrary{babel,patterns,decorations.pathreplacing}
\usepackage{dsfont} 

\usepackage{setspace} 

\usepackage{booktabs} 
\usepackage{paralist}

\deffootnote{1.8em}{0em}{\textsuperscript{\thefootnotemark}\ }

\makeatletter
\g@addto@macro\bfseries{\boldmath}
\makeatother 

\newtheoremstyle{mythm}
{\topsep}   
{\topsep}   
{\itshape}  
{0pt}       
{\bfseries} 
{.}         
{5pt plus 1pt minus 1pt} 
{}          

\newtheoremstyle{mydefi}
{\topsep}   
{3ex}   
{\normalfont}  
{0pt}       
{\bfseries} 
{.}         
{5pt plus 1pt minus 1pt} 
{}          

\makeatletter
\newenvironment{myproof}[1][\proofname]{\par
	\pushQED{\qed}%
	\normalfont \topsep6\p@\@plus6\p@\relax
	\trivlist
	\item[\hskip\labelsep
	\itshape
	#1\@addpunct{.}]\ignorespaces
}{%
	\popQED\endtrivlist\@endpefalse\bigskip
}
\makeatother

\theoremstyle{mythm}
\newtheorem{satz}{Satz}[section]
\newtheorem{thm}[satz]{Theorem}
\newtheorem{prop}[satz]{Proposition}
\newtheorem{lem}[satz]{Lemma}

\newtheorem{cor}[satz]{Corollary}

\newtheorem{rem}[satz]{Remark}

\theoremstyle{mydefi}
\newtheorem{defi}[satz]{Definition}

\newtheorem{ex}[satz]{Example}

\DeclareMathOperator{\Aut}{Aut}
\DeclareMathOperator{\SL}{SL}
\DeclareMathOperator{\GL}{GL}
\DeclareMathOperator{\PSL}{PSL}

\DeclareMathOperator{\PGammaL}{P\Gamma L}
\DeclareMathOperator{\GammaL}{\Gamma L}
\DeclareMathOperator{\N}{N}
\DeclareMathOperator{\PG}{PG}

\newcommand{\inv}[1]{#1^{-1}}
\newcommand{\mat}[4]{\left( \begin{smallmatrix}#1&#2\\#3&#4\end{smallmatrix} \right)}

\usepackage[backend=biber,hyperref]{biblatex}
\usepackage[colorlinks=true,linkcolor=.,citecolor=.,urlcolor=blue,breaklinks=true,%
pdftitle={Automorphisms of (Affine) SL(2,q)-Unitals},%
pdfauthor={Verena Möhler}%
]{hyperref}

\title{Automorphisms of (Affine) \texorpdfstring{$\SL(2,q)$-Unitals}{SL(2,q)-Unitals}}
\author{Verena Möhler}
\date{December 17, 2020}

\begin{document}

\maketitle

\begin{abstract}
	\noindent $\SL(2,q)$-unitals are unitals of order $q$ admitting a regular action of $\SL(2,q)$ on the 
	complement of some block. They can be obtained from affine $\SL(2,q)$-unitals via parallelisms. We compute a sharp upper bound for automorphism groups of affine $\SL(2,q)$-unitals and show that exactly two parallelisms are fixed by all automorphisms. 
	In $\SL(2,q)$-unitals obtained as closures of affine $\SL(2,q)$-unitals via those two parallelisms, we show that there is one block fixed under the full automorphism group.\bigskip
	
	\noindent \textbf{2020 MSC:}
	51A10, 
	05E18 
	\smallskip
	
	\noindent \textbf{Keywords:} design, unital, affine unital, automorphism, parallelism
\end{abstract}

Most of the results in the present paper have been obtained in the author's Ph.\,D.\ thesis \cite{diss}, where detailed arguments can be found for some statements that
we leave to the reader here.

\section{Preliminaries}

A \textbf{unital of order $n$} is a $2$-$(n^3+1,n+1,1)$ design, i.\,e.\ an incidence structure with $n^3+1$ points, $n+1$ points on each block and unique joining blocks for any two points. We also consider \emph{affine unitals}, which arise from unitals by removing one block (and all the points on it) and can be completed to unitals via a parallelism on the short blocks. We give an axiomatic description:

\begin{defi}
	Let $n\in\mathbb N$, $n\geq2$. An incidence structure $\mathbb U=(\mathcal P,\mathcal B,\emph I)$ is called an \textbf{affine unital of order \boldmath $n$} if:
	\begin{itemize}
		\item[(AU1)] There are $n^3-n$ points.
		\item[(AU2)] Each block is incident with either $n$ or $n+1$ points. The blocks incident with $n$~points will be called \textbf{short blocks} and the blocks incident with $n+1$ points will be called \textbf{long blocks}.
		\item[(AU3)] Each point is incident with $n^2$ blocks.
		\item[(AU4)] For any two points there is exactly one block incident with both of them.
		\item[(AU5)] There exists a \textbf{parallelism} on the short blocks, meaning a partition of the set of all short blocks into $n+1$ parallel classes of size $n^2-1$ such that the blocks of each parallel class are pairwise non-intersecting.
	\end{itemize}	
\end{defi}

The existence of a parallelism as in (AU5) must explicitly be required (see \cite[Example 3.10]{diss}). An affine unital $\mathbb U$ of order $n$ with parallelism $\pi$ can be completed to a unital $\mathbb U^\pi$ of order $n$ as follows: For each parallel class, add a new point that is incident with each short block of that class. Then add a single new block $[\infty]^\pi$, incident with the $n+1$ new points (see \cite[Proposition 3.9]{diss}). We call $\mathbb U^\pi$ the \textbf{$\pi$-closure} of $\mathbb U$.\bigskip

Note that though we must require the existence of a parallelism in the definition of an affine unital, this parallelism need not be unique. It is therefore not convenient to require that isomorphisms of affine unitals respect certain parallelisms and we will only ask them to be isomorphisms of the incidence structures. We call two parallelisms $\pi$ and $\pi'$ of an affine Unital $\mathbb U$ \textbf{equivalent} if there is an automorphism of $\mathbb U$ which maps $\pi$ to~$\pi'$.\bigskip

Given an affine unital $\mathbb U$ with parallelisms $\pi$ and $\pi'$, the closures $\mathbb U^\pi$ and $\mathbb U^{\pi'}$ are isomorphic with $[\infty]^\pi \mapsto [\infty]^{\pi'}$ exactly if $\pi$ and $\pi'$ are equivalent (see \cite[Proposition 3.12]{diss}). Indeed, two closures $\mathbb U^\pi$ and $\mathbb U^{\pi'}$ being isomorphic does not imply the parallelisms $\pi$ and $\pi'$ being equivalent in $\mathbb U$, as is shown with the following

\begin{ex}
Let $\mathbb U$ be the affine unital of order $3$ (i.\,e.\ having 24 points, 30 long blocks and 32 short blocks) indicated in figure \ref{zweiparas}.
\begin{figure}
\centering
\begin{tikzpicture}[scale=0.23]
\foreach \x in {0,...,62} \draw (\x,0) -- (\x,-24) (\x,-24.5) -- (\x,-28.5) (\x,-29) -- (\x,-33);
\foreach \y in {0,...,-24} \draw (0,\y) -- (62,\y);
\foreach \y in {-24.5,...,-28.5} \draw (0,\y) -- (62,\y);
\foreach \y in {-29,...,-33} \draw (0,\y) -- (62,\y);

\foreach \y/\x in {{ 1, 4, 5, 6 }/0, { 1, 7, 8, 9 }/1, { 1, 12, 13, 14 }/2, { 1, 15, 16, 17 }/3, { 1, 22, 23, 24 }/4, { 2, 4, 13, 15 }/
	5, { 2, 8, 10, 23 }/6, { 2, 11, 18, 21 }/7, { 2, 12, 19, 22 }/8, { 2, 16, 20, 24 }/9, { 3, 5, 14, 17 }/10,
	{ 3, 6, 9, 18 }/11, { 3, 7, 12, 21 }/12, { 3, 8, 11, 16 }/13, { 3, 10, 19, 20 }/14, { 4, 8, 18, 24 }/15,
	{ 4, 9, 12, 23 }/16, { 4, 14, 16, 21 }/17, { 5, 10, 16, 18 }/18, { 5, 13, 20, 22 }/19, { 5, 15, 19, 21 }/20,
	{ 6, 7, 19, 24 }/21, { 6, 8, 17, 20 }/22, { 6, 10, 15, 22 }/23, { 7, 11, 17, 22 }/24, { 7, 14, 20, 23 }/25,
	{ 9, 10, 13, 21 }/26, { 9, 11, 14, 19 }/27, { 11, 12, 15, 24 }/28, { 13, 17, 18, 23 }/29,
	{ 1, 2, 3 }/30, { 1, 10, 11 }/51, { 1, 18, 19 }/59, { 1, 20, 21 }/38, { 2, 5, 7 }/54, { 2, 6, 14 }/46, { 2, 9, 17 }/
	39, { 3, 4, 22 }/40, { 3, 13, 24 }/55, { 3, 15, 23 }/47, { 4, 7, 10 }/31, { 4, 11, 20 }/60, { 4, 17, 19 }/52,
	{ 5, 8, 12 }/43, { 5, 9, 24 }/48, { 5, 11, 23 }/35, { 6, 11, 13 }/44, { 6, 12, 16 }/36, { 6, 21, 23 }/56,
	{ 7, 13, 16 }/49, { 7, 15, 18 }/41, { 8, 13, 19 }/37, { 8, 14, 15 }/57, { 8, 21, 22 }/50, { 9, 15, 20 }/32,
	{ 9, 16, 22 }/58, { 10, 12, 17 }/61, { 10, 14, 24 }/42, { 12, 18, 20 }/53, { 14, 18, 22 }/33, { 16, 19, 23 }/45,
	{ 17, 21, 24 }/34
}
{\foreach \z in \y \draw[fill=black] (\x,-\z) rectangle (\x+1,-\z+1);}
\foreach \x in {30,...,37} \draw[pattern=north east lines] (\x,-25.5) rectangle (\x+1,-24.5);
\foreach \x in {38,...,45} \draw[pattern=north east lines] (\x,-26.5) rectangle (\x+1,-25.5);
\foreach \x in {46,...,53} \draw[pattern=north east lines] (\x,-27.5) rectangle (\x+1,-26.5);
\foreach \x in {54,...,61} \draw[pattern=north east lines] (\x,-28.5) rectangle (\x+1,-27.5);
\foreach \x in {30,...,34,43,44,45} \draw[fill=black!20] (\x,-30) rectangle (\x+1,-29);
\foreach \x in {30,...,34,43,44,45} \draw[pattern=dots] (\x,-30) rectangle (\x+1,-29);
\foreach \x in {38,...,42,35,36,37} \draw[fill=black!20] (\x,-31) rectangle (\x+1,-30);
\foreach \x in {38,...,42,35,36,37} \draw[pattern=dots] (\x,-31) rectangle (\x+1,-30);
\foreach \x in {46,...,50,59,60,61} \draw[fill=black!20] (\x,-32) rectangle (\x+1,-31);
\foreach \x in {46,...,50,59,60,61} \draw[pattern=dots] (\x,-32) rectangle (\x+1,-31);
\foreach \x in {54,...,58,51,52,53} \draw[fill=black!20] (\x,-33) rectangle (\x+1,-32);
\foreach \x in {54,...,58,51,52,53} \draw[pattern=dots] (\x,-33) rectangle (\x+1,-32);

\foreach \x in {1,...,24} \draw[font=\tiny] (-1.5,-\x+0.5) node {\x};
\foreach \x in {1,3,...,61} \draw[font=\tiny] (\x-0.5,0.5) node {\x};
\foreach \x in {2,4,...,62} \draw[font=\tiny] (\x-0.5,1.5) node {\x};
\foreach \x in {25,...,28} \draw[font=\tiny] (-1.5,-\x) node {\x} (-1.5,-\x-4.5) node {\x};

\draw[decorate,decoration={brace,amplitude=3pt}] (-2.5,-28.5)--(-2.5,-24.5) node[midway,left]{$\pi$};
\draw[decorate,decoration={brace,amplitude=3pt}] (-2.5,-33)--(-2.5,-29) node[midway,left]{$\pi'$};
\end{tikzpicture}
\caption{An affine unital with two non-equivalent parallelisms, where the closures are isomorphic}
\label{zweiparas}
\end{figure}
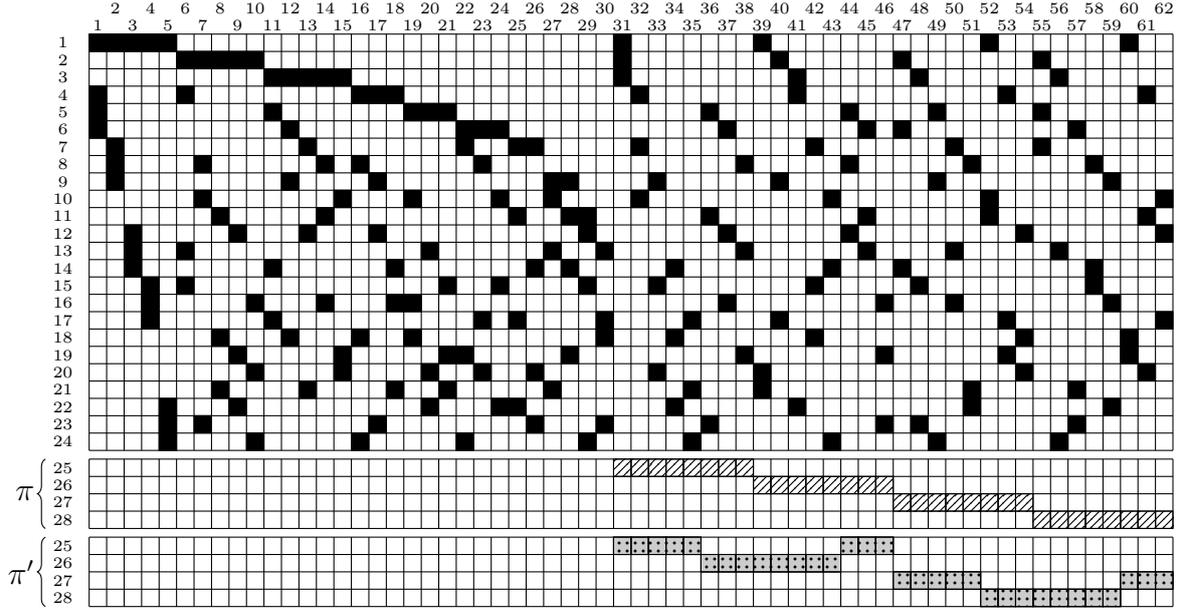
We consider two parallelisms $\pi$ and $\pi'$ of $\mathbb U$, in the figure indicated via the incidences of the $4$ additional points $25,\ldots,28$. Those two parallelisms are not equivalent in $\mathbb U$ -- in fact, $\Aut(\mathbb U)$ is trivial -- but the corresponding closures $\mathbb U^\pi$ and $\mathbb U^{\pi'}$ are isomorphic via the isomorphism
\[ (1,16,23,10)(2,11,15,19)(4,9,13,14)(5,22,18,24)(6,27,20,25,8,26,17,28)(12,21) \text{.}\]

\end{ex}

\section{(Affine) \texorpdfstring{$\SL(2,q)$-Unitals}{SL(2,q)-Unitals}}

From now on let $p$ be a prime and $q\coloneqq p^e$ a $p$-power. We are interested in a special kind of affine unitals, namely affine $\SL(2,q)$-unitals. The construction of those affine unitals is due to Grundhöfer, Stroppel and Van Maldeghem \cite{slu}. They consider translations of unitals, i.\,e.\ automorphisms fixing each block through a given point (the so-called center). Of special interest are unitals of order $q$ where two points are centers of translation groups of order $q$. In the classical (Hermitian) unital of order $q$, any two such translation groups generate a group isomorphic to $\SL(2,q)$; see \cite[Main Theorem]{moufang} for further possibilities. The construction of (affine) $\SL(2,q)$-unitals is motivated by this action of $\SL(2,q)$ on the classical unital.\bigskip

Let $S\leq \SL(2,q)$ be a subgroup of order $q+1$ and let $T\leq \SL(2,q)$ be a Sylow $p$-subgroup. Recall that $T$ has order $q$ (and thus trivial intersection with $S$), that any two conjugates $T^h\coloneqq\inv hTh$, $h\in \SL(2,q)$, have trivial intersection unless they coincide and that there are $q+1$ conjugates of $T$.

Consider a collection $\mathcal D$ of subsets of $\SL(2,q)$ such that each set $D\in \mathcal D$ contains \mbox{$\mathds{1}\coloneqq\mat 1001$}, that $\# D=q+1$ for each $D\in \mathcal D$, and the following properties hold:

\begin{itemize}
	\item[(Q)] For each $D\in \mathcal D$, the map
	\[(D\times D)\smallsetminus \{(x,x) \mid x\in D\} \to \SL(2,q)\text{,}\quad (x,y)\mapsto x\inv y\text{,}\]
	is injective, i.\,e.\ the set $D^*\coloneqq\{x\inv y \mid x,y\in D,\ x\neq y\}$
	contains $q(q+1)$ elements.
	\item[(P)] The system consisting of $S\smallsetminus\{\mathds{1}\}$, all conjugates of $T\smallsetminus\{\mathds{1}\}$ and all sets $D^*$ with $D\in \mathcal D$ forms a partition of $\SL(2,q)\smallsetminus\{\mathds{1}\}$.
\end{itemize}

Set
\begin{align*}
\mathcal P \coloneqq &\ \SL(2,q)\text{,}\\
\mathcal B \coloneqq &\ \{Sg \mid g\in \SL(2,q)\} \cup \{T^h g \mid h,g\in \SL(2,q)\} \cup \{Dg \mid D\in\mathcal D, g\in\SL(2,q)\}
\end{align*}
and let the incidence relation $I\subseteq \mathcal P\times\mathcal B$ be containment.\bigskip

\goodbreak

Then we call the incidence structure $\mathbb U_{S,\mathcal D}\coloneqq (\mathcal P,\mathcal B,I)$ an \textbf{\boldmath affine $\SL(2,q)$-unital}. Each affine $\SL(2,q)$-unital is indeed an affine unital of order $q$, see \cite[Prop. 3.15]{diss}. If $\pi$ is a parallelism on the short blocks of an affine $\SL(2,q)$-unital $\mathbb U_{S,\mathcal D}$, we call the $\pi$-closure~$\mathbb U_{S,\mathcal D}^\pi$ an \textbf{$\SL(2,q)$-($\pi$-)unital}.\bigskip

In any affine $\SL(2,q)$-unital, the set of short blocks is the set of all right cosets of the Sylow $p$-subgroups, i.\,e.\
\[ \{Tg \mid T\in\mathfrak P, g\in\SL(2,q)\}\text{,}\]
where $\mathfrak P$ denotes the set of the $q+1$ Sylow $p$-subgroups of $\SL(2,q)$.
Note that each right coset $Tg$ is a left coset $gT^g$ of a conjugate of $T$. A parallelism as in (AU5) means a partition of the set of short blocks into $q+1$ sets of $q^2-1$ pairwise non-intersecting cosets. For each prime power $q$, there are hence two obvious parallelisms, namely partitioning the set of short blocks into the sets of \emph{right} cosets or into the sets of \emph{left} cosets of the Sylow $p$-subgroups. We name those two parallelisms ``flat'' and ``natural'', respectively, and denote them by the corresponding musical signs
\[ \flat \coloneqq \{\{Tg\mid g\in\SL(2,q)\}\mid T\in\mathfrak P\}\quad \text{and}\quad \natural \coloneqq \{\{gT\mid g\in\SL(2,q)\}\mid T\in\mathfrak P\}\text{.} \]

\goodbreak

\begin{ex}\leavevmode
	\begin{enumerate}[(a)]
		\item For each prime power $q$ we may choose $S=C$ to be cyclic and $\mathcal H$ a set of blocks through $\mathds{1}$ such that $\mathbb U_{C,\mathcal H}$ is isomorphic to the affine part of the classical unital and the closure $\mathbb U_{C,\mathcal H}^\natural$ is isomorphic to the classical unital. We call $\mathbb U_{C,\mathcal H}$ the \textbf{classical affine $\SL(2,q)$-unital} and $\mathbb U_{C,\mathcal H}^\natural$ the \textbf{classical $\SL(2,q)$-unital}. See \cite[Example 3.1]{slu} or \cite[Section 3.2.2]{diss} for details.
		
		In \cite[Example 5.5]{slu}, Grundhöfer, Stroppel and Van Maldeghem state that the unital $\mathbb U_{C,\mathcal H}^\flat$ coincides with a unital described by Grüning, that can be embedded in a Hall plane and in its dual, see \cite{gruening}.
		\item Several non-classical affine $\SL(2,q)$-unitals are known, namely one of order $4$, described in \cite{slu}, and three of order $8$, described in
		\cite[Section 3]{threeaff}.
		\item There are also other parallelisms than $\flat$ and $\natural$, see \cite[Sections 2.1-2.3]{paratrans}. The construction of $\SL(2,q)$-unitals inspired Nagy and Mez{\H o}fi to a method of constructing new incidence structures from old ones by removing a block and attaching it again in a different way, see \cite{paramod}. By computer search, they found plenty parallelisms for affine unitals of orders $3$ and $4$, respectively. Among the parallelisms of affine $\SL(2,q)$-unitals, we show in Theorem \ref{flatundnat} that $\flat$ and $\natural$ are the only ones being preserved by all possible automorphisms of affine $\SL(2,q)$-unitals.
	\end{enumerate}
\end{ex}

\section{Automorphisms}
Concerning automorphisms of (affine) $\SL(2,q)$-unitals, we note first that on any affine $\SL(2,q)$-unital $\mathbb U_{S,\mathcal D}$, the group $\SL(2,q)$ acts as group of automorphisms by multiplication from the right. Each automorphism of $\SL(2,q)$ induces a permutation on the point set of $\mathbb U_{S,\mathcal D}$, while it need not leave the block set invariant. Since any automorphism of affine unitals maps short blocks to short blocks, we will first consider the incidence structure~$\mathfrak S$ given by the points and short blocks of any affine $\SL(2,q)$-unital. Recall that the short blocks are the same in every affine $\SL(2,q)$-unital, independent of the choice of $S$ and the set $\mathcal D$.

\subsection{Automorphisms of the Geometry of Short Blocks}

We use an embedding of $\mathfrak S$ into the classical generalized quadrangle $\mathrm{Q}(4,q)$. This quadrangle is the polar space given by the equation
\[ x_1x_3 + x_2x_4 +x_5^2 = 0\]
in homogeneous coordinates for the projective space $\PG(4,q)$. The intersection of $\mathrm{Q}(4,q)$ with the hyperplane $x_5=0$ of $\PG(4,q)$ is a geometric hyperplane $H$ of the generalized quadrangle, isomorphic to the classical generalized quadrangle $\mathrm Q(3,q)$, which is the polar space given by the equation
\[ x_1x_3 + x_2x_4 = 0 \] 
in homogeneous coordinates for the projective space $\PG(3,q)$. The geometry of short blocks $\mathfrak S$ is isomorphic to the complement $\mathrm Q(4,q)\smallsetminus H$, which follows basically from \cite[\nopp 10.7.8]{paynethas} and \cite[Theorem 1.1]{kthas} (see \cite[Section 4.1]{diss} for details). In any thick classical generalized quadrangle, the automorphisms of a hyperplane complement are exactly the automorphisms of the quadrangle stabilizing the hyperplane (see \cite[Lemma 2.3]{pasini}). Thus,
\[ \Aut(\mathfrak S) \cong \Aut(\mathrm Q(4,q))_H \text{.} \]
Compute $\Aut(\mathrm Q(4,q))_H\cong \mathrm O(f) \rtimes \Aut(\mathbb F_q)$,
where $f$ is the quadratic form
\[ f\colon \mathbb F_q^4\to \mathbb F_q\text{, } x=(x_1,\ldots,x_4)\mapsto x_1x_3 + x_2x_4\text{,} \]
and $\mathrm O(f)= \{A\in\GL(4,q)\mid f(xA)=f(x)\ \forall x\in\mathbb F_q^4\}$ is the corresponding orthogonal group (again, see \cite[Section 4.1]{diss} for details).\bigskip

We thus know the isomorphism type (and, in particular, the order) of $\Aut(\mathfrak S)$ and may now study the action of $\Aut(\mathfrak S)$ on the incidence structure $\mathfrak S$. We already know that $\SL(2,q)$ acts via right multiplication as a group of automorphisms on $\mathfrak S$. Applying automorphisms of $\SL(2,q)$ to the point set of $\mathfrak S$ also induces automorphisms of $\mathfrak S$ and so does inversion, since the set of short blocks of any affine $\SL(2,q)$-unital is given by the set of \emph{all} right cosets of the Sylow $p$-subgroups of $\SL(2,q)$, and every left coset of a Sylow $p$-subgroup $T$ is a right coset of a conjugate of $T$. (Recall $\inv{(Tg)}=\inv g T = T^g\inv g$.)\bigskip

Let $R \coloneqq \{\rho_h \mid h\in\SL(2,q)\} \leq \Aut(\mathfrak S)$, where $\rho_h \in R$ acts on $\mathfrak S$ by right multiplication with $h\in \SL(2,q)$. Let further $\mathfrak A\leq \Aut(\mathfrak S)$ denote the subgroup given by automorphisms of $\SL(2,q)$ and $I \leq\Aut(\mathfrak S)$ the cyclic subgroup of order $2$ given by inversion. For each $a\in\GL(2,q)$, conjugation by $a$ induces an automorphism $\gamma_a\in\mathfrak A$.

\begin{lem}
	$\Aut(\mathfrak S) = (\mathfrak A\times I)\cdot R$.
\end{lem}

\begin{myproof}
	From the above, the product $\mathfrak A\cdot I\cdot R$ is a subset  of $\Aut(\mathfrak S)$. The product $\mathfrak A\cdot I$ is direct (and in particular a subgroup of $\Aut(\mathfrak S)$), since $\mathfrak A$ and $I$ commute and have trivial intersection. The group $R$ has trivial intersection with $\mathfrak A\times I$, since $\mathfrak A$ and $I$ fix the point~$\mathds{1}$, while $R$ acts regularly on the point set. Hence we may compute
	\begin{align*}
	\#((\mathfrak A\times I)\cdot R) &= \#\mathfrak A \cdot \#I \cdot \#R\\
	&= \#\PGammaL(2,q) \cdot 2 \cdot \#\SL(2,q) \\
	&= e(q-1)q(q+1) \cdot 2 \cdot (q-1)q(q+1) \\
	&= 2e(q-1)^2 q^2 (q+1)^2\text{.}
	\end{align*}
	We have already noted that $\Aut(\mathfrak S)$ is isomorphic to $\mathrm O (f) \ltimes \Aut(\mathbb F_q)$. But $\#\mathrm O (f) = 2(q-1)^2q^2(q+1)^2$, see \cite[141]{taylor} or \cite[Theorems 6.21 and 7.23 with $\nu=2$ and $\delta=0$]{wan}, and the lemma follows.	
\end{myproof}

\begin{lem}\label{arsemidir}
	$\mathfrak A$ normalizes $R$ in the full automorphism group of $\mathfrak S$.
\end{lem}

\begin{myproof}
	Let $\alpha\in \mathfrak A$, $\rho_h\in R$ and let $x\in \SL(2,q)$ be a point of $\mathfrak S$. Then
	\[x\cdot(\inv\alpha \rho_h \alpha)= ((x\cdot\inv\alpha)h)\cdot \alpha = x(h\cdot\alpha) = x\cdot \rho_{h\cdot\alpha}\text{.}\qedhere\]
\end{myproof}

Since $\mathfrak A$ and $R$ have trivial intersection, we thus know that the product of $\mathfrak A$ and $R$ is semidirect.

\subsection{Automorphisms of Affine \texorpdfstring{$\SL(2,q)$-Unitals}{SL(2,q)-Unitals}}

We use our knowledge of the full automorphism group of $\mathfrak S$ to compute automorphism groups of (affine) $\SL(2,q)$-unitals.

\begin{thm}\label{isom}  Let $q\geq 3$ and let $\mathbb U_{S,\mathcal D}$ and $\mathbb U_{S',\mathcal D'}$ be affine $\SL(2,q)$-unitals.
	\begin{enumerate}[\em (a)]
		\item \label{isoma} Let $\psi\colon \mathbb U_{S,\mathcal D} \to \mathbb U_{S',\mathcal D'}$ be an isomorphism. Then $\psi = \alpha\rho_h$ with $\rho_h\in R$ and $\alpha\in \mathfrak A$ such that $S\cdot \alpha=S'$.
		\item \label{isomb} $\Aut(\mathbb U_{S,\mathcal D})\leq \mathfrak A_S\ltimes R$. \qed
	\end{enumerate}
\end{thm}

\begin{myproof}
	\begin{enumerate}[(a)]
		\item Each isomorphism of affine $\SL(2,q)$-unitals maps short blocks to short blocks and is hence an automorphism of $\mathfrak S$, since the incidence structure of the short blocks is the same in every affine $\SL(2,q)$-unital. We thus know $\psi \in \Aut(\mathfrak S) = (\mathfrak A\times I)\cdot R$. Let $\psi = \alpha \iota \rho_h$ with $\iota\in I$ inversion or identity. Then 
		\[ S\cdot \psi = S\cdot(\alpha \iota \rho_h) = (S\cdot\alpha)h\text{,} \]
		since $S$ is a group and $\alpha$ an automorphism of $\SL(2,q)$. Since right multiplication by $\SL(2,q)$ induces automorphisms in every affine $\SL(2,q)$-unital, $S\cdot\alpha$ is a block through $\mathds{1}$ in $\mathbb U_{S',\mathcal D'}$ and is a subgroup of $\SL(2,q)$ of order $q+1$. The only block with these properties is $S'$ and hence $S\cdot\alpha= S'$. Assume $\iota$ to be inversion and choose $g\in\SL(2,q)$ with $g\cdot\alpha\not\in \N(S')$. Note that such an element $g$ exists, since $q\geq 3$. Then
		\[
		(Sg)\cdot \psi = (Sg)\cdot(\alpha \iota \rho_h)=
		((S\cdot\alpha)(g\cdot\alpha))^{-1}h = (S')^{g\cdot\alpha}(g\cdot\alpha)^{-1}h\text{.}	
		\]
		$Sg$ is a block of $\mathbb U_{S,\mathcal D}$ and hence $(S')^{g\cdot\alpha}$ is a block of $\mathbb U_{S',\mathcal D'}$ with $\mathds{1}$ in $(S')^{g\cdot\alpha}$ and $(S')^{g\cdot\alpha}\leq \SL(2,q)$ is a subgroup of order $q+1$. But $(S')^{g\cdot\alpha}\neq S'$ by the choice of~$g$, a contradiction.
		\item Directly from Lemma \ref{arsemidir} and (\ref{isoma}). \qedhere
	\end{enumerate}
\end{myproof}

\begin{rem}
	The classical affine $\SL(2,q)$-unital $\mathbb U_{C,\mathcal H}$ admits the whole group $\mathfrak A_C\ltimes R$ as automorphism group \emph{(}see \emph{\cite[Proposition 4.6]{diss})}. Hence, $\mathfrak A_S\ltimes R$ is a sharp upper bound for the automorphism group of any affine $\SL(2,q)$-unital $\mathbb U_{S,\mathcal D}$ of order $q\geq 3$.
\end{rem}

\begin{ex}
	It is easily seen that there is only one isomorphism type of affine unitals of order $2$.
	Labelling the points with elements of $\SL(2,2)$, we get the classical affine $\SL(2,2)$-unital $\mathbb U_{C,\mathcal H}$ with $C=\langle \mat 0111\rangle$ and $\mathcal H$ the empty set:
	
	
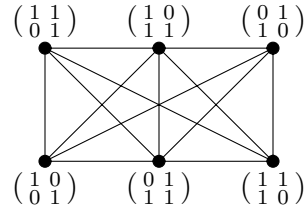
\begin{center}
		\begin{tikzpicture}[scale=1.5]
		\draw[fill=black] (0,0) circle (1.5pt) node[below]{$\mat 1001$}-- (1,0) circle (1.5pt) node[below]{$\mat 0111$}-- (2,0) circle (1.5pt) node[below]{$\mat 1110$}-- (2,1) circle (1.5pt) node[above]{$\mat 0110$}--(1,1) circle (1.5pt) node[above]{$\mat 1011$}--(0,1) circle (1.5pt) node[above]{$\mat 1101$};
		\draw (0,1)--(0,0)--(1,1) --(2,0)--(0,1)--(1,0)--(2,1)--(0,0) (1,0)--(1,1);
		\end{tikzpicture}
	\end{center}

	Obviously, inversion induces an automorphism of this affine unital (interchanging the points $\mat 0111$ and $\mat 1110$) and we get
	$\Aut(\mathbb U_{C,\mathcal H}) = (\mathfrak A_C \times I)\cdot R \gneq \mathfrak A_C \ltimes R$
	for $q=2$.
	
\end{ex}\bigskip

\subsection{Automorphisms of \texorpdfstring{$\SL(2,q)$-Unitals}{SL(2,q)-Unitals} with parallelism \texorpdfstring{$\flat$ or $\natural$}{flat or natural}}

We now take a closer look at $\SL(2,q)$-unitals with parallelism $\flat$ or $\natural$. Recall the definitions
\[ \flat \coloneqq \{\{Tg\mid g\in\SL(2,q)\}\mid T\in\mathfrak P\}\quad \text{and}\quad \natural \coloneqq \{\{gT\mid g\in\SL(2,q)\}\mid T\in\mathfrak P\}\text{,} \]
where $\mathfrak P$ denotes the set of all Sylow $p$-subgroups of $\SL(2,q)$. Considering the action of~$R$ (i.\,e.\ right multiplication by $\SL(2,q)$), we see that both $\flat$ and $\natural$ are invariant under this action. On $\flat$, right multiplication by $\SL(2,q)$ fixes each parallel class, while on $\natural$ right multiplication with $h\in\SL(2,q)$ maps any left coset $gT$ to the left coset $ghT^h$ and hence the group $R$ acts on the parallel classes of $\natural$ via conjugation on the Sylow $p$-subgroups. Since the group $\mathfrak A$ acts via automorphisms of $\SL(2,q)$, both parallelisms $\flat$ and $\natural$ are fixed under the action of $\mathfrak A$. Note that inversion interchanges $\flat$ with $\natural$.\bigskip

Since the stabilizer of the block $[\infty]=[\infty]^\pi$ in the full automorphism group of the $\pi$-closure of any affine unital $\mathbb U$ equals the group of automorphisms of $\mathbb U$ fixing the parallelism $\pi$, we get the following

\begin{cor}
	Let $q\geq 3$, $\pi\in\{\flat,\natural\}$ and let $\mathbb U_{S,\mathcal D}$ be an affine $\SL(2,q)$-unital. Then \[\Aut(\mathbb U_{S,\mathcal D}^\pi)_{[\infty]} = \Aut(\mathbb U_{S,\mathcal D})\] and in particular $\Aut(\mathbb U_{C,\mathcal H}^\pi)_{[\infty]} = \mathfrak A_C\ltimes R$ for $q\geq 3$. \qed
\end{cor}

The two parallelisms $\flat$ and $\natural$ are indeed the only parallelisms in any affine $\SL(2,q)$-unital that are preserved under the action of $R$:

\begin{thm}\label{flatundnat}
	Let $\mathbb U_{S,\mathcal D}$ be an affine $\SL(2,q)$-unital with parallelism $\pi$ such that right multiplication by $\SL(2,q)$ preserves $\pi$. Then $\pi\in\{\flat,\natural\}$.
\end{thm}

\begin{myproof}
	Let $T\in\mathfrak P$ and let $[T]$ denote the parallel class of $\pi$ containing the short block~$T$. Since right multiplication with $t\in T$ fixes the block $T$, it fixes also the parallel class~$[T]$. We distinguish two cases.
	
	Assume first that $[T]$ contains exactly the $q^2-1$ left cosets of $T$. Then for each $h\in\SL(2,q)$, we have 
	\[ [T]\cdot h=\{gT\mid g\in\SL(2,q)\}\cdot h = \{gTh \mid g\in\SL(2,q)\} = \{ghT^h \mid g\in\SL(2,q)\} \text{,} \]
	which equals the set of left cosets of $T^h$ and must be a parallel class of $\pi$. Hence, each parallel class of $\pi$ is the complete set of left cosets of some Sylow $p$-subgroup, meaning $\pi=\natural$.
	
	Assume now that $[T]$ contains at least one short block $T^hg=gT^{hg}$ which is not a left coset of $T$, i.\,e.\ $T^{hg}\neq T$. Then right multiplication with $T$ as well as right multiplication with~$T^{hg}$ fixes a block in $[T]$ and hence also the parallel class $[T]$. Since any two different Sylow $p$-subgroups generate $\SL(2,q)$, we have $\langle T\cup T^{hg}\rangle =\SL(2,q)$. Hence, right multiplication with $\SL(2,q)$ stabilizes $[T]$ and $[T]$ equals the set of right cosets of $T$. The same reasoning works for each parallel class and thus $\pi = \flat$.
\end{myproof}

\begin{defi}\label{deftrans}
	A \textbf{translation} with \textbf{center} $c$ of a unital $\mathbb U$ is an automorphism of~$\mathbb U$ that fixes the point $c$ and each block through $c$. The group of all translations with center~$c$ will be denoted by $G_{[c]}$.
	
	We call $c$ a \textbf{translation center} if $G_{[c]}$ acts transitively on the set of points different from~$c$ on any block through $c$.
\end{defi}

\begin{rem}\label{semireg}
	Let $\mathbb U$ be a unital of order $q$ and $c$ a point of $\mathbb U$. Then the group $G_{[c]}$ of all translations of $\mathbb U$ with center $c$ acts semiregularly on the points of $\mathbb U$ different from $c$ \emph{(}see \emph{\cite[Theorem 1.3]{alltrans}}\emph{)}. Hence, $\#G_{[c]}\leq q$ and $c$ is a translation center exactly if $\#G_{[c]}=q$.
\end{rem}

In any $\SL(2,q)$-$\pi$-unital $\mathbb U_{S,\mathcal D}^\pi$, we label the points on the block $[\infty]$ with the Sylow $p$-subgroups in such a way that each (affine short) block $T\in\mathfrak P$ through~$\mathds{1}$ is incident with the point $T\in[\infty]$.

\begin{lem}\label{sylow_trans}
	Let $\mathbb U_{S,\mathcal D}^\natural$ be an $\SL(2,q)$-$\natural$-unital. For each $T\in\mathfrak P$, the point $T\in[\infty]$ is a translation center with $G_{[T]}=R_T\coloneqq\{\rho_t\mid t\in T\}$.
\end{lem}

\begin{myproof} 
	For each $T\in\mathfrak P$, the set of blocks through the point $T$ in $\mathbb U_{S,\mathcal D}^\natural$ equals $\{[\infty]\}\cup\{gT\mid g\in\SL(2,q)\}$. Hence, $R_T$ obviously is a group of translations of $\mathbb U_{S,\mathcal D}^\natural$ with center $T$. Since $\# R_T=q$, the statement follows (recall Remark \ref{semireg}).
\end{myproof}

We use this statement on translations and a theorem of Grundhöfer, Stroppel and Van Maldeghem \cite{alltrans} to show that the block $[\infty]$ is fixed by every automorphism in any non-classical $\SL(2,q)$-$\natural$-unital.

\begin{prop} \label{naturalinfstab}
	Let $\mathbb U_{S,\mathcal D}^\natural$ be a non-classical $\SL(2,q)$-$\natural$-unital. Then \[\Aut(\mathbb U_{S,\mathcal D}^\natural) = \Aut(\mathbb U_{S,\mathcal D}^\natural)_{[\infty]}\text{.} \]
\end{prop}

\begin{myproof} 
	Assume that there is an automorphism of $\mathbb U_{S,\mathcal D}^\natural$ not fixing $[\infty]$. Since all the points on~$[\infty]$ are translation centers (see Lemma \ref{sylow_trans}), there are thus three non-collinear translation centers of $\mathbb U_{S,\mathcal D}^\natural$. Then $\mathbb U_{S,\mathcal D}^\natural$ is the classical unital, as is shown in \cite{alltrans}.
\end{myproof}

Grüning showed that in each Grüning unital $\mathbb U_{C,\mathcal H}^\flat$, the block $[\infty]$ is fixed by every automorphism (see \cite[Lemma 5.5]{gruening}). We will extend this statement to every $\SL(2,q)$-$\flat$-unital $\mathbb U_{S,\mathcal D}^\flat$ of order $q\geq 3$, independent of the group $S$ and the set $\mathcal D$. For order $2$, there is only one isomorphism type of unitals, represented by the classical unital. Hence, any unital of order $2$ admits a $2$-transitive automorphism group and no block is fixed by the full automorphism group.\bigskip

Our proof of Theorem \ref{flatinfstab} below uses further knowledge about the group $S\leq \SL(2,q)$ of order $q+1$ and its stabilizer in $\Aut(\SL(2,q))$. The following statements are verified by using known results about $\PSL(2,q)$, see \cite[Section 2.2]{diss} for details.

\begin{prop}[\cite{diss}, Proposition 2.5]\label{gruppes}
	Let $S\leq\SL(2,q)$ be a subgroup of order $q+1$. Then:
	\begin{enumerate}[\em (a)]
		\item For $q\not\equiv 3 \mod 4$, the group $S$ is cyclic.
		\item For $q\equiv 3 \mod 4$, there are the following possibilities:
		\begin{enumerate}[\em (i)]
			\item $S$ is cyclic.
			\item $S$ is a generalized quaternion group and the quotient $S/\{\pm \mathds{1}\}$ in $\PSL(2,q)$ is a dihedral group of order $\frac 12(q+1)$.
			\item The quotient $S/\{\pm \mathds{1}\}$ in $\PSL(2,q)$ is either isomorphic to $A_4$ with $q=23$, or isomorphic to $S_4$ with $q=47$. \qed
		\end{enumerate}		
	\end{enumerate}
\end{prop}

\goodbreak

\begin{thm}[\cite{diss}, Theorem 2.11] \label{autsls}\leavevmode
	\begin{enumerate}[\em (a)]
		\item For $C\leq \SL(2,q)$ cyclic of order $q+1$, we have \[\Aut(\SL(2,q))_C \cong C_{q+1} \rtimes C_{2e}\text{.} \]
		\item For $q\equiv 3 \mod 4$, $q>7$ and $S\leq \SL(2,q)$ generalized quaternion of order $q+1$, we~have \[\#\Aut(\SL(2,q))_S = e(q+1)\]
		and $\Aut(\SL(2,q))_S$ is conjugate to a subgroup of $\Aut(\SL(2,q))_C$ with $C$ cyclic of order $q+1$.
		\item For $q\in\{23,47\}$ and $S\leq\SL(2,q)$ of type \emph{(b)(iii)} in \emph{Proposition \ref{gruppes}} or $q=7$ and $S\leq\SL(2,7)$ a quaternion group, we have
		\[\pushQED{\qed} \Aut(\SL(2,q))_S \cong S_4\text{.} \qedhere \popQED\]  
	\end{enumerate}
\end{thm}

\begin{cor} \label{autslscor}
	Let $S\leq \SL(2,q)$ be a subgroup of order $q+1$. Then the stabilizer $\Aut(\SL(2,q))_S$ is solvable.\qed
\end{cor}

\begin{defi}
	Let $G$ be a group. Using the \textbf{commutator series} $(G^{(n)})_{n\in\mathbb N_0}$, we define the \textbf{stable commutator} of $G$ by $G^{(\omega)}\coloneqq\bigcap_{n\in\mathbb N_0} G^{(n)}$.
\end{defi}

\begin{thm} \label{flatinfstab}
	Let $\mathbb U_{S,\mathcal D}^\flat$ be an $\SL(2,q)$-$\flat$-unital of order $q\geq 3$. Then every automorphism of $\mathbb U_{S,\mathcal D}^\flat$ fixes the block $[\infty]$.
\end{thm}

\begin{myproof}
	We treat the small cases $q=3$ and $q=4$ first:	
	For $q=3$, there is only one isomorphism type of affine $\SL(2,q)$-unitals (see \cite[Theorem 3.3]{slu}), namely the classical affine $\SL(2,q)$-unital. Hence, its $\flat$-closure is the Grüning unital and the theorem holds by \cite[Lemma 5.5]{gruening}.
	
	For $q=4$, an exhaustive computer search using GAP showed that there are only two isomorphism types of affine $\SL(2,q)$-unitals (see \cite[Section 6.1.1]{diss}) and none of their $\flat$-closures admits automorphisms which move the block $[\infty]$.\bigskip
	
	Recall that right multiplication by $\SL(2,q)$ fixes each parallel class of the parallelism~$\flat$ and that hence the group $R\leq\Aut(\mathbb U_{S,\mathcal D}^\flat)_{[\infty]}$ fixes each point on $[\infty]$. Recall further that we label the points on $[\infty]$ with the Sylow $p$-subgroups in such a way that each (affine short) block $T\in\mathfrak P$ through~$\mathds{1}$ is incident with the point $T\in[\infty]$.
	
	Assume that there exists an automorphism $\alpha \in \Aut(\mathbb U_{S,\mathcal D}^\flat)$ with $[\infty]\cdot\alpha\neq [\infty]$. Since $R$ acts transitively on the affine part of $\mathbb U_{S,\mathcal D}^\flat$, assume without restriction $\mathds{1}\in [\infty]\cdot\alpha$. We~distinguish two cases.
	
	\begin{description}
		\item[Case 1:] The block $[\infty]\cdot\alpha$ intersects $[\infty]$, i.\,e.\ $[\infty]\cdot\alpha = T$ for $T\leq \SL(2,q)$ a Sylow $p$-subgroup.
	\end{description}
	
	Assume first that $\alpha$ does not fix $T$ as a point. Then $[\infty]\cdot\inv\alpha$ meets $[\infty]$ in the point $T\cdot\inv\alpha\neq T$. Choose $h\in\SL(2,q)$ such that $([\infty]\cdot\inv\alpha)h\neq [\infty]\cdot\inv\alpha$. Then the automorphism $\inv\alpha\rho_h\alpha$ moves $[\infty]$ and we have $T\cdot(\inv\alpha\rho_h\alpha) = T$ as a point, since $T\cdot\inv\alpha \in [\infty]$.
	
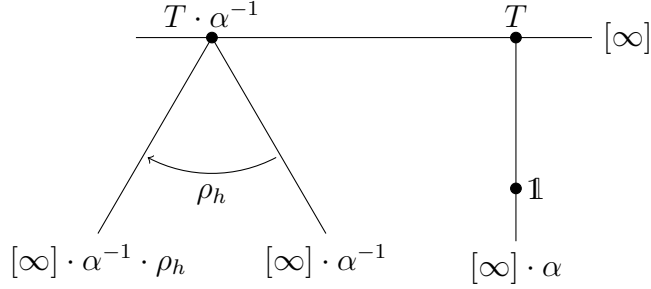
	\begin{center}
		\begin{tikzpicture}
		\draw[fill=black] (0,0) coordinate (a) circle (2pt) node[above] {$T\cdot\inv\alpha$} (4,0) coordinate (b) circle(2pt) node[above] {$T$} (4,-2) circle(2pt) node[right]{$\mathds{1}$};
		\draw (-1,0)--(5,0) node[right] {$[\infty]$};
		\draw (b)--+(0,-2.7) node[below]{$[\infty]\cdot\alpha$} (a)--(240:3) node[below]{$[\infty]\cdot\inv\alpha\cdot\rho_h$} (a)--(300:3) node[below]{$[\infty]\cdot\inv\alpha$};
		\draw[<-] (242:1.8) arc(242:298:1.8) node[midway, below]{$\rho_h$};	
		\end{tikzpicture}
	\end{center}
	
	We may thus assume without restriction that $\alpha$ fixes $T$ as a point. The group
	\[R_T\coloneqq \{\rho_t \in R \mid t\in T\}\leq R\]
	acts regularly on the affine points of the block $T$ and trivially on $[\infty]$. Hence, the group of automorphisms $\alpha R_T\inv\alpha$ acts regularly on $[\infty]\smallsetminus\{T\}$ and trivially on $[\infty]\cdot\inv\alpha$. In particular, an affine point is fixed and we have
	\[q=\#\alpha R_T\inv\alpha \mid \#\Aut(\mathbb U_{S,\mathcal D}^\flat)_{[\infty],\mathds{1}}\mid \# \Aut(\SL(2,q))_S\text{.} \]
	
	According to Theorem \ref{autsls}, this implies $q\mid 2e(q+1)$ or $q\mid e(q+1)$, since $24=\# S_4$ is not divided by $7$, $23$ or $47$, respectively. Now $\gcd(q,q+1)=1$ and $q=p^e>e$ is obtained by an easy induction on $e$. Hence, $q=p^e\mid 2e$. Since $p^e > e$, this implies $p=2$ and $2^{e-1}\mid e$. Again, induction shows that $2^{e-1}>e$ if $e\geq 3$. It remains $q\in \{2,4\}$, but we are only interested in $q\geq 3$ and have already excluded the case $q=4$.
	
	\begin{description}
		\item[Case 2:] The block $[\infty]\cdot\beta$ does not meet $[\infty]$ for any automorphism $\beta\in \Aut(\mathbb U_{S,\mathcal D}^\flat)$ with $[\infty]\cdot\beta\neq[\infty]$.
	\end{description}

	Recall the assumption $\mathds{1} \in [\infty]\cdot\alpha$. Since $[\infty]\cdot\alpha$ does not intersect $[\infty]$, we know $[\infty]\cdot\alpha \neq T$ for any Sylow $p$-subgroup $T\leq \SL(2,q)$.
	
	We show that no two blocks in the orbit $[\infty]\cdot\Aut(\mathbb U_{S,\mathcal D}^\flat)$ may intersect in one point: Let~$\vartheta\in\Aut(\mathbb U_{S,\mathcal D}^\flat)$ such that $[\infty]\cdot\alpha\neq [\infty]\cdot\alpha\vartheta$ and assume without restriction that the two blocks $[\infty]\cdot\alpha$ and $[\infty]\cdot\alpha\vartheta$ intersect in the point $\mathds{1}$. Then the automorphism $\alpha\vartheta\inv\alpha$ moves the block $[\infty]$. Let $x\coloneqq \mathds{1}\cdot\inv\vartheta \in[\infty]\cdot\alpha$. Note that, other than indicated in the picture, $x$ need not be be different from $\mathds{1}$. Then $(x\cdot\inv\alpha)\in[\infty]$ and
	\[ (x\cdot\inv\alpha)\cdot (\alpha\vartheta\inv\alpha) = \mathds{1}\cdot\inv \alpha\in [\infty]\text{.} \]
	Hence, the automorphism $\alpha\vartheta\inv\alpha$ moves the block $[\infty]$ to a block intersecting $[\infty]$, a~contradiction to the requirement.
	
	
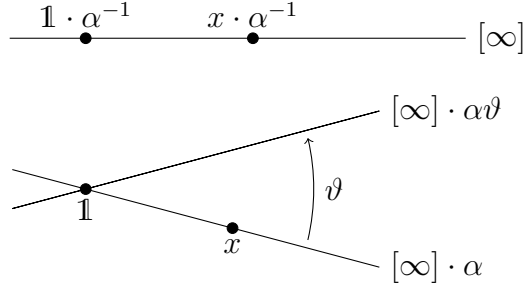
\begin{center}
		\begin{tikzpicture}
		\draw[fill=black] (0,0) coordinate (a) circle (2pt) node[below] {$\mathds{1}$} (-15:2) coordinate (b) circle(2pt) node[below] {$x$} (0,2) circle(2pt) node[above]{$\mathds{1}\cdot\inv\alpha$} (2.2,2) circle(2pt) node[above]{$x\cdot\inv\alpha$};
		\draw (-1,2)--(5,2) node[right] {$[\infty]$};
		\draw (-15:-1)--(-15:4) node[right]{$[\infty]\cdot\alpha$} (15:-1)--(15:4) node[right]{$[\infty]\cdot\alpha\vartheta$} (15:-1)--(15:4);
		\draw[->] (-13:3) arc(-13:13:3) node[midway, right]{$\vartheta$};	
		\end{tikzpicture}
	\end{center}
	
	For any block $D\in\mathcal D$ through $\mathds{1}$ and $\mathds{1}\neq d\in D$, the two blocks $D$ and $D\inv d$ intersect in the point $\mathds{1}$. Since right multiplication by $\inv d$ is an automorphism of $\mathbb U_{S,\mathcal D}^\flat$, the block $[\infty]\cdot\alpha$ may thus be none of the blocks $D\in\mathcal D$. Hence, $[\infty]\cdot\alpha=S$. Since $R\leq\Aut(\mathbb U_{S,\mathcal D}^\flat)$ and every affine point is contained in a right coset of $S$, we have
	\[ \mathcal O_{[\infty]}\coloneqq[\infty]\cdot\Aut(\mathbb U_{S,\mathcal D}^\flat) = \{Sg\mid g\in\SL(2,q)\}\cup \{[\infty]\}\text{.} \]
	
	\begin{center}
		\begin{tikzpicture}
		\foreach \k in {0,1,2,4,5,6} \draw (0,-0.5*\k)--(5,-0.5*\k);
		\draw (2.5,-1.4) node {$\vdots$};
		\draw (5,0) node[right]{$[\infty]$} (5,-3) node[right]{$S$};
		\draw [decorate,decoration={brace,amplitude=6pt}] (5.6,-0.4) -- (5.6,-3.1) node[midway,xshift=1.5cm] {$q^2-q$ blocks};
		\draw[->] (-0.2,0) .. controls (-1,-1.5) .. (-0.2,-3) node[midway,left]{$\alpha$};
		\draw[fill=black] (1,-3) circle (2pt) node[below]{$\mathds{1}$};
		\end{tikzpicture}
	\end{center}
	
	The group of automorphisms $R$ fixes $[\infty]$ and acts transitively on the $q^2-q$ right cosets of $S$ and hence $\Aut(\mathbb U_{S,\mathcal D}^\flat)$ acts $2$-transitively on $\mathcal O_{[\infty]}$. Since \[\#\mathcal O_{[\infty]} = q^2-q+1\text{,} \] the group $\Aut(\mathbb U_{S,\mathcal D}^\flat)$ must be a $2$-transitive group on $q^2-q+1$ elements and contain \mbox{$R\cong \SL(2,q)$}. We use the classification of finite $2$-transitive groups to reach a contradiction.
	
	There are several sources for the classification of finite $2$-transitive groups, of which the most convenient one seems to be two lists by Cameron, see \cite[Tables 7.3 and 7.4]{cameron}. In~the majority of cases, the degree of the action already yields a contradiction, why we only copied the information needed in Tables \ref{camlistaff} and \ref{camlistalmsimp}. In both tables, $f$ is a prime power. In the list of affine $2$-transitive groups, $G_0$ denotes the stabilizer of one point in the $2$-transitive action. In the list of almost simple $2$-transitive groups, $N$ denotes the minimal normal subgroup of the $2$-transitive group $G$. Indeed, there is one action of degree $9^2=3^4$ missing in Cameron's list of affine $2$-transitive groups, where $\SL(2,5)$ is a normal subgroup of $G_0$, see e.\,g.\ Liebeck's proof of Hering's classification of affine $2$-transitive groups \cite[Appendix 1]{liebeck}. However, we are only interested in the degree in this case and the missing action will not be relevant for our argument. \bigskip
	
	\begin{table}
		\centering
		\caption{Degrees of affine 2-transitive groups}
		\label{camlistaff}
		\begin{tabular}{lll}\toprule
			\textbf{Degree}	& \boldmath{$G_0$}				& \textbf{Condition}\\\midrule
			$f^d$	& $\SL(d,f)\leq G_0\leq \GammaL(d,f)$	&					\\
			$f^{2d}$&										&					\\
			$f^6$	&										&					\\
			$f$		&						&$f=5^2,7^2,11^2,23^2,3^4,11^2,19^2,$\\
			&						&$\phantom{f=}\quad 29^2,59^2,2^4,2^6,3^6$\\\bottomrule
		\end{tabular}
	\end{table}
	
	\begin{table}
		\centering
		\caption{Degrees of almost simple 2-transitive groups}
		\label{camlistalmsimp}
		\begin{tabular}{lll}\toprule
			\textbf{Degree}	& \boldmath{$N$}			& \textbf{Condition}\\\midrule
			$n$				& $A_n$						& $n\geq 5$				\\
			$(f^d-1)/(f-1)$	& $\PSL(d,f)$				& $d\geq 2$, $(d,f)\neq(2,2),(2,3)$\\
			$2^{2d-1}\pm2^{d-1}$&							& $d\geq 3$					\\
			$f^3+1$		&								& 					\\
			$f^2+1$		&								&					\\
			$f$	&							&$f=11,12,15,22,23,24,28,176,276$\\\bottomrule		
		\end{tabular}
	\end{table}
	
	Let $\Gamma_{[\infty]}\coloneqq\Aut(\mathbb U_{S,\mathcal D}^\flat)_{[\infty]}$ denote the stabilizer of $[\infty]\in\mathcal O_{[\infty]}$ in our required $2$-transitive action of $\Aut(\mathbb U_{S,\mathcal D}^\flat)$ on $\mathcal O_{[\infty]}$. Then 
	\[\SL(2,q) \cong R \leq \Gamma_{[\infty]} \leq \mathfrak A_S\ltimes R\cong \Aut(\SL(2,q))_S\ltimes \SL(2,q)\text{.} \]
	
	In particular, we know that the stable commutator $\Gamma_{[\infty]}^{(\omega)}$ equals $\SL(2,q)$ for $q\geq 4$, since $\Aut(\SL(2,q))_S$ is solvable for each choice of $S$, see Corollary \ref{autslscor}. In the first affine $2$-transitive action in Table \ref{camlistaff}, the stable commutator $G_0^{(\omega)}$ equals the stable commutator of $\SL(d,f)$, which is trivial or $\SL(d,f)$. Since all isomorphisms between finite (projective) special linear groups are known (see e.\,g.\ \cite[Satz 6.14]{huppert}), we conclude $d=2$ and $f=q$. But this is not possible because of our required degree $q^2-q+1\neq q^2$ and we may thus exclude the first affine $2$-transitive action.
	
	\goodbreak
	
	Every other affine $2$-transitive action is excluded by the observation that their degree is always a square but our degree $q^2-q+1$ is never a square, since
	\[ (q-1)^2 < q^2-q+1 < q^2\text{.} \]
	
	Now we look at the list of almost simple $2$-transitive groups. The point stabilizer of the action of $A_n$ is $A_{n-1}$, which is far bigger than $\Gamma_{[\infty]}$ for $n=q^2-q+1$.
	
	The second entry in Table \ref{camlistalmsimp} is the action on the points or hyperplanes of a projective space. For $d\geq 3$ and $(d-1,f)\not\in \{(2,2),(2,3)\}$, the stable commutator of the point stabilizer of this action is the special affine group 
	\[\mathrm{ASL}(d-1,f)= \SL(d-1,f)\ltimes \mathbb F_f^{d-1}\text{.} \]
	Hence, $\mathrm{ASL}(d-1,f)\cong \SL(2,q)$ (recall $q\geq 4$). But $\mathrm{ASL}(d-1,f)$ contains the normal subgroup $\mathbb F_f^{d-1}$ of order $f^{d-1}> 4$, while the highest possible order of a normal subgroup of $\SL(2,q)$ is $2$ for $q\geq 4$.
	For $(d-1,f)\in \{(2,2),(2,3)\}$ the degree of the action would be $7$ resp.\ $13$ and hence $q=3$ or $q=4$, which we already excluded in the beginning. It~remains the case $d=2$. Then the degree is $f+1$ and hence $f=q^2-q=q(q-1)$. But $\gcd(q,q-1)=1$ shows that the product $q(q-1)$ is not a prime power. Hence we may exclude the second entry in the list of almost simple $2$-transitive actions.\bigskip
	
	For the remaining cases, considering the degree will suffice. The degrees $2^{2d-1}\pm2^{d-1}$ ($d\geq 3$) and $12,22,24,28,176,276$ are excluded by the observation that $q^2-q+1$ is odd. Then we just showed that $q^2-q$ cannot be a prime power $f^d$. Finally, it is easy to see that $q^2-q=q(q-1)\not \in \{10,14,22\}$ for any prime power $q$.\bigskip
	
	Thus, $\Aut(\mathbb U_{S,\mathcal D}^\flat)$ does not act $2$-transitively on $\mathcal O_{[\infty]}$ and our proof is complete.	
\end{myproof}\bigskip

\textbf{Acknowledgment.} The author wishes to warmly thank her thesis advisor Markus J. Stroppel for his highly valuable support in each phase of this research.
	
	\printbibliography
\end{document}